\documentclass{amsart}

\usepackage{amsfonts}
\usepackage{amsmath}

\theoremstyle{theorem}
\newtheorem{lem}{Lemma}
\newtheorem{prop}[lem]{Proposition}
\newtheorem{cor}[lem]{Corollary}
\newtheorem{thm}[lem]{Theorem}

\theoremstyle{definition}

\newtheorem{defn}[lem]{Definition}

\title{Uniformizable and Realcompact Bornological Universes}
\author{Tom Vroegrijk}
\date{}

\begin{document}

\maketitle

\begin{abstract}
Bornological universes were introduced by Hu in \cite{Hu1} and obtained renewed interest in recent articles on convergence in hyperspaces and function spaces and optimization theory. In \cite{Hu1} and \cite{Hu2} Hu gives a necessary and sufficient condition for which a bornological universe is metrizable. In this article we will give a characterization of uniformizable bornological universes. Furthermore, a construction on bornological universes that the author used in \cite{Vroegrijk2} to find the bornological dual of function spaces endowed with the bounded-open topology will be used to define realcompactness for bornological universes. We will also give various characterizations of this new concept.
\end{abstract}

% Keywords: bornology; uniform space; totally bounded; realcompactness.
% MSC: 54E15, 54E35, 54C25

\section{Introduction}

A \emph{bornological universe} $(X,\mathcal{B})$ consists of a topological space $X$ and a bornology $\mathcal{B}$ on the underlying set of $X$. Bornological universes play a key role in recent publications on convergence structures on hyperspaces \cite{Beer5, Beer6, Beer7, Lechicki}, topologies on function spaces \cite{Beer4} and optimization theory \cite{Beer2, Beer3}. In \cite{Hu1} Hu defines a bornological universe $(X,\mathcal{B})$ to be metrizable if there is a metric $d$ on $X$ such that the topology defined by $d$ is the original topology on $X$ and the elements of $\mathcal{B}$ are exactly the sets that are bounded for this metric. In this article the following characterization of metrizable bornological universes is given: 
\begin{thm}
A bornological universe $(X,\mathcal{B})$ is metrizable iff the following conditions are satisfied:
\begin{enumerate}
\item $X$ is metrizable,
\item $\mathcal{B}$ has a countable base,
\item for each $B_1\in\mathcal{B}$ there is a $B_2\in\mathcal{B}$ such that $\overline{B_1}\subseteq \overset{\circ}{B_2}$
\end{enumerate}
\end{thm}
\noindent Keeping this result of Hu in mind we can ask ourselves when a bornological universe is uniformizable. The first problem that we encounter is how to define uniformizability for bornological universes. Each uniform space has an underlying topological space, but there are (at least) two natural bornologies that we can associate with a uniformity: the bornology of sets that are bounded in the sense of Bourbaki and the bornology of totally bounded sets. This gives us two possible ways to define uniformizable bornological universes. We will see, however, that both are equivalent. One of the various characterizations of uniformizability for bornological universes that we will give is being isomorphic to a subspace of a product of real lines. This property was also studied by Beer in \cite{Beer1}, but his definition of a product bornology is different from the one that we will use here.\\

In \cite{Vroegrijk2} the author uses generalized bornologies (see \cite{Vroegrijk1}) to extend the duality between locally convex topological vector spaces and bornological vector spaces to vector spaces with a topology defined by extended quasinorms. One of the objects that is encountered in this article is the space $\mathcal{C}(X)$ of continuous, real valued maps on a Tychonoff space $X$ endowed with the bounded-open topology for some bornology $\mathcal{B}$ on $X$. To describe the bornological dual of this object a realcompact extension $\upsilon_{\mathcal{B}}(X)$ of $X$ is introduced. The way this space $\upsilon_{\mathcal{B}}(X)$ is defined for a bornological universe $(X,\mathcal{B})$ is similar to the definition of the Hewitt realcompactification for topological spaces. We will use this to define realcompactness for bornological universes.

\section{Uniformizable bornological universes}

\begin{defn}
A \emph{bornology} on a set $X$ is a set $\mathcal{B}\subseteq 2^X$ that satisfies the following conditions:
\begin{enumerate}
\item[B1] $\{x\}\in\mathcal{B}$ for each $x\in X$,
\item[B2] if $B\in\mathcal{B}$ and $A\subseteq B$, then $A\in\mathcal{B}$,
\item[B3] $A\cup B\in\mathcal{B}$ whenever $A,B\in\mathcal{B}$.
\end{enumerate}
The elements of a bornology are called \emph{bounded sets} and a map that preserves boundedness is called a \emph{bounded map}. A topological space $X$ endowed with a bornology $\mathcal{B}$ will be called a bornological universe.
\end{defn}

It is easily verified that the category of bornological universes and bounded, continuous maps is a topological construct. A source $(f_i:(X,\mathcal{B})\to (X_i,\mathcal{B}_i))_{i\in I}$ in this category is initial if $X$ is endowed with the topology that is initial for the source $(f_i:X\to X_i)_{i\in I}$ and a set $B\subseteq X$ is bounded iff $f_i(B)$ is bounded for each $i\in I$. This means that if $(X,\mathcal{B})$ is a bornological universe and $Y\subseteq X$, the subspace structure on $Y$ consists of the subspace topology on $Y$ and the bornology of all subsets of $Y$ that are bounded in $X$. The product structure on a product of a family of bornological universes $(X_i,\mathcal{B}_i)_{i\in I}$ on the other hand is defined as the product topology and the bornology of sets $B\subseteq \prod_{i\in I}X_i$ for which each projection $\pi_i(B)$ is bounded.

\begin{defn}
A subset $B$ of a uniform space $(X,\mathcal{U})$ is called \emph{bounded in the sense of Bourbaki} (see \cite{Bourbaki}) if for each entourage $U\in\mathcal{U}$ we can find an $n\in\mathbb{N}$ and a finite set $K\subseteq X$ such that $B\subseteq U^n(K)$.
\end{defn}

For a locally convex topological vector space the classical notion of boundedness coincides with boundedness in the sense of Bourbaki with respect to the canonical uniformity. In a metric space $(X,d)$ we have that each set that is bounded for the metric $d$ is bounded in the sense of Bourbaki for the underlying uniformity, but the converse is in general not true.\\

If $(X,\mathcal{U})$ is a uniform space, then the sets that are bounded in the sense of Bourbaki form a bornology on $X$. A second bornology that we can associate with the uniformity $\mathcal{U}$ is the bornology of totally bounded subsets. A subset $B$ of a uniform space is \emph{totally bounded} if for each entourage $U$ there is a finite partition $(B_i)_{i=1}^n$ of $B$ such that $B_i\times B_i\subseteq U$ for each $i\in I$. This is equivalent with saying that for each entourage $U\in\mathcal{U}$ there is a finite set $K\subseteq X$ such that $B\subseteq U(K)$, so obviously each totally bounded set is bounded in the sense of Bourbaki. The proof of the following theorem can be found in \cite{Hejcman}.

\begin{thm}
Let $(X,\mathcal{U})$ be a uniform space. The following statements are equivalent:
\begin{enumerate}
\item $A$ is bounded in the sense of Bourbaki,
\item $A$ is bounded for each uniformly continuous pseudometric $d$ on $X$,
\item each real valued, uniformly continuous map is bounded on $A$.
\end{enumerate}
\end{thm}

The bornology of sets that are bounded in the sense of Bourbaki and the bornology of totally bounded sets are two natural bornologies associated with a uniform space. At first sight it therefore looks like we can define two notions of uniformizability for topological bornological spaces. We will see that they are actually equal.

\begin{defn}
We will call a bornological universe $(X,\mathcal{B})$ \emph{$bb$-uniformizable} if there is a uniform space $(X,\mathcal{U})$ for which $X$ is the underlying topological space and $\mathcal{B}$ is the set of subsets of $X$ that are bounded in the sense of Bourbaki for $\mathcal{U}$. In the case that $\mathcal{B}$ is equal to the set of totally bounded subsets of the uniform space $(X,\mathcal{U})$ we will call $(X,\mathcal{B})$ \emph{$tb$-uniformizable}.
\end{defn}

Let $(X,\mathcal{U})$  be a uniform space and take $\mathcal{V}\subseteq \mathcal{U}$. If each entourage in $\mathcal{U}$ contains a finite intersection of elements of $\mathcal{V}$, then $\mathcal{V}$ is called a \emph{base} for $\mathcal{U}$. For a set $\mathcal{C}$ of real valued maps on a set $X$ we will denote the uniform space that is initial for the source $\mathcal{C}$ as $(X,\mathcal{C})$. Define $U_f^{\epsilon}$ as the set that contains all $(x,y)\in X\times X$ for which $|f(x)-f(y)|<\epsilon$. The entourages $U_f^{\epsilon}$, form a base for the uniformity of $(X,\mathcal{C})$.

\begin{lem}\label{totallyboundedbase}
Let $\mathcal{V}$ be a base for a uniformity $\mathcal{U}$. If for each $V\in\mathcal{V}$ we can find a finite partition $(B_i)_{i=0}^n$ of $B$ such that $B_i\times B_i\subseteq V$ for all $0\leq i\leq n$, then $B$ is totally bounded.
\begin{proof}
Each entourage $U\in\mathcal{U}$ contains a set $V_0\cap \ldots \cap V_m$ where each $V_k$ is in $\mathcal{V}$. For all $0\leq k\leq m$ we can find a finite partition $\mathcal{P}_k$ of $B$ such that $B_k\times B_k\subseteq V_k$ for each $B_k\in\mathcal{P}_k$. If we define $\mathcal{P}$ as $\{B_0\cap \ldots \cap B_m|\forall 0\leq k\leq m:B_k\in \mathcal{P}_k\}$, then $\mathcal{P}$ is a finite partition of $B$ and each of each of its elements $P$ satisfies $P\times P\subseteq U$.
\end{proof}
\end{lem}

\begin{lem}\label{boundedinitial}
For a set $\mathcal{C}$ of realvalued maps on a set $X$, the following statements are equivalent:
\begin{enumerate}
\item $f(B)$ is bounded for all $f\in\mathcal{C}$,
\item $B$ is totally bounded in $(X,\mathcal{C})$,
\item $B$ is bounded in the sense of Bourbaki in $(X,\mathcal{C})$.
\end{enumerate}
\begin{proof}
If we assume that the first statement is true and we take $\epsilon>0$ and $f\in\mathcal{C}$, then we can find a finite partition $(C_i)_{i=0}^n$ of $f(B)$ such that each element of this partition has diameter smaller than $\epsilon$. The sets $B_i$, defined as $B\cap f^{-1}(C_i)$, form a finite partition of $B$ and for each $0\leq i\leq n$ holds $B_i\times B_i\subseteq U_f^{\epsilon}$. Applying lemma \ref{totallyboundedbase} gives us that $B$ is totally bounded in $(X,\mathcal{C})$.

We already saw that the second statement implies the third. Now suppose that $B$ is bounded in the sense of Bourbaki. This yields that it is bounded for each uniformly continuous metric on $(X,\mathcal{C})$.  Since $d_f$, where $d_f(x,y)$ is defined as $|f(x)-f(y)|$, is a uniformly continuous metric for each $f\in\mathcal{C}$, we have that $f(B)$ is bounded.
\end{proof}
\end{lem}

\begin{defn}\label{satborn}
We will call a bornology $\mathcal{B}$ on a topological space \emph{saturated} if it contains each $B$ that satisfies the following condition:
if $(G_n)_n$ is a sequence of open sets such that each bounded set is contained in some $G_n$ and $\overline{G_n}\subseteq G_{n+1}$ for all $n\in\mathbb{N}$, then there is an $m\in\mathbb{N}$ for which $B\subseteq G_m$.
\end{defn}

From here on the set of all continuous, real valued maps on $X$ that are bounded on all elements of $\mathcal{B}$ will be denoted as $\mathcal{C}_{\mathcal{B}}(X)$.

\begin{prop}
The following statements are equivalent:
\begin{enumerate}
\item $\mathcal{C}_{\mathcal{B}}(X)$ is an initial source,
\item $(X,\mathcal{B})$ is $bb$-uniformizable,
\item $X$ is completely regular and $\mathcal{B}$ is saturated.
\end{enumerate}
\begin{proof}
If the first statement is true, then we have by definition that the underlying topological space of $(X,\mathcal{C}_{\mathcal{B}}(X))$ is the original topology in $X$. We know from lemma \ref{boundedinitial} that the set of subsets of $X$ that are bounded in the sense of Bourbaki in $(X,\mathcal{C}_{\mathcal{B}}(X))$ is exactly $\mathcal{B}$.

Now suppose that $(X,\mathcal{B})$ is $bb$-uniformizable. We automatically obtain that $X$ is completely regular. If $B$ is an unbounded subset of $X$, then we can find a uniformly continuous metric $d$ on $X$ for which $B$ is unbounded. Take an arbitrary $x_0\in X$ and define $G_n$ as the open ball with center $x_0$ and radius $n+1$. This is an increasing sequence of open sets that satisfies the conditions stated in definition \ref{satborn}. Furthermore, we know that $B$ is not contained in any of the sets $G_n$. Hence we have that $\mathcal{B}$ is saturated.

For a completely regular space $X$, the initial topology for the source that consists of all continuous maps into $[0,1]$ is the original topology. Clearly, these maps are all in $\mathcal{C}_{\mathcal{B}}(X)$. If we now assume that the third statement is true and that $B$ is an unbounded set, then we can find a sequence $(G_n)_n$ such that each bounded set is contained in some $G_n$, $\overline{G_n}\subseteq G_{n+1}$ for all $n\in\mathbb{N}$ and no $G_n$ contains $B$. Choose a sequence $(x_n)_n$ in $B$ such that $x_n\in B\setminus \overline{G_n}$. Let $\phi_n$ be a map into $[0,n]$ that vanishes on $\overline{G_n}$ and attains the value $n$ in $x_n$. Define the map $\phi$ as $\sum_n\phi_n$. Since $\phi$ is equal to $\sum_{k=0}^n\phi_k$ on each $G_n$ we obtain that $\phi$ is a continuous map. By definition, $\phi$ is bounded on all elements of $\mathcal{B}$ and unbounded on $B$.
\end{proof}
\end{prop}

A subset of a topological space is called \emph{relatively pseudocompact} if it is mapped to a bounded subset of $\mathbb{R}$ by all real valued, continuous maps on $X$. The previous proposition yields that each relatively pseudocompact subset of a uniformizable bornological universe is bounded.

\begin{prop}
An object $(X,\mathcal{B})$ is $bb$-uniformizable iff it is $tb$-uni\-formizable.
\begin{proof}
If $(X,\mathcal{B})$ is $bb$-uniformizable, then $\mathcal{C}_{\mathcal{B}}(X)$ is initial. From lemma \ref{boundedinitial} we obtain that $(X,\mathcal{C}_{\mathcal{B}}(X))$ is a uniform space with underlying topological space $X$ for which the totally bounded subsets are exactly the sets in $\mathcal{B}$.

Now let $(X,\mathcal{B})$ be a $tb$-uniformizable object and $B$ an unbounded subset of $X$. We can find  a uniformly continuous metric $d$ on $X$ and a sequence $(x_n)_n$ in $B$ such that for distinct numbers $n$ and $m$ it holds that $d(x_n,x_m)\geq 1$. Define $K_n$ as the set $\{x_m|m\geq n\}$ and $G_n$ as the set of all $x\in X$ that lie at a $d$-distance strictly greater than $1/(n+3)$ from $K_n$. This is an increasing sequence of open sets and $\overline{G_n}\subseteq G_{n+1}$ for all $n\in\mathbb{N}$. This sequence also satisfies that each bounded set is contained in some $G_n$. If this were not the case, then there would be a bounded set that contains a countable subset with the property that all of its elements, for the metric $d$, lie at a distance greater than $1/3$ of each other. This would imply that this bounded set is not totally bounded. Since $B$ is not contained in any of the sets $G_n$ we obtain that $\mathcal{B}$ is saturated.
\end{proof}
\end{prop}

From here on a $bb$-uniformizable bornological universe will be simply called \emph{uniformizable}. We will say that a bornological universe $(X,\mathcal{B})$ is \emph{Hausdorff} if  $X$ is Hausdorff.

\begin{prop}\label{reallinessubspace}
The following statements are equivalent:
\begin{enumerate}
\item $(X,\mathcal{B})$ is Hausdorff and uniformizable,
\item $(X,\mathcal{B})$ is isomorphic to a subspace of a product of real lines.
\end{enumerate}
\begin{proof}
If $(X,\mathcal{B})$ is Hausdorff uniformizable, then the source $\mathcal{C}_{\mathcal{B}}(X)$ is initial and separating. This yields that the map from $X$ to $\mathbb{R}^{\mathcal{C}_{\mathcal{B}}(X)}$ that sends an element $x$ to $(f(x))_{f\in\mathcal{C}_{\mathcal{B}}(X)}$, is in fact an embedding.

To prove the converse, we need to show that each subspace of a product of real lines is Hausdorff uniformizable. We can endow each subset $X$ of a product $\mathbb{R}^{\alpha}$ of real lines with the uniformity that it inherits from the product uniformity on $\mathbb{R}^{\alpha}$. This uniformity is the initial one for the source that consists of all projection maps to $\mathbb{R}$. We  know that the underlying topology of this uniformity is the original topology on $X$. Lemma \ref{boundedinitial} grants us that a subset of $X$ is bounded in the sense of Bourbaki for this uniformity iff each of its projections is bounded.
\end{proof}
\end{prop}

An article by Beer (see \cite{Beer1}) also contains a necessary and sufficient condition for a bornological universe to be embeddable into a product of real lines. The product bornology, however, is in this article defined as the bornology generated by all sets for which at least one projection is bounded, while we use the bornology of sets for which all projections are bounded.

\section{Realcompact bornological universes}

A subset $B$ of a locally convex topological vector space $E$ is called \emph{bounded} if each continuous seminorm is bounded on $B$. If $E$ satisfies the condition that each seminorm that is bounded on all bounded sets is automatically continuous, then $E$ is called \emph{bornological}. The bornological objects form a concretely coreflective subcategory of the category of locally convex topological vector spaces.\\

Let $\mathcal{B}$ be a bornology on a Tychonoff space $X$ that consist only of relatively pseudocompact sets. The set of all real valued, continuous maps on $X$ endowed with the bounded-open topology is a locally convex topological vector space. In \cite{Schmets} Schmets gives a characterization of the bornological coreflection of this topological vector space using the Hewitt realcompactification $\upsilon(X)$ of $X$.\\

If the elements of $\mathcal{B}$ are no longer supposed to be relatively pseudocompact, then the bounded-open topology is no longer a vector topology. In \cite{Vroegrijk2}, however, an extension of the classical duality between topological and bornological vector spaces is given that allows us to define the bornological coreflection of this space. In that same article a characterization of this bornological coreflection is given using a realcompact space $\upsilon_{\mathcal{B}}(X)$ that contains $X$ as a dense subspace. We will use these spaces to define realcompactness for bornological universes and give various characterizations of this new concept.\\

Let $(X,\mathcal{B})$ be a bornological universe where $X$ is a Tychonoff space. The space $\upsilon_{\mathcal{B}}(X)$ is defined as $\{x\in\beta(X)|\forall f\in\mathcal{C}_{\mathcal{B}}(X): f(x)\not=\infty\}$. Here $f^{\beta}$ is the unique map from $\beta(X)$ to the one-point compactification of $\mathbb{R}$ that satisfies the condition that its restriction to $X$ is equal to $f$.

\begin{prop}
$\upsilon_{\mathcal{B}}(X)$ is a realcompact Tychonoff space that contains each $B\in\mathcal{B}$ as a relatively compact subspace, i.e. the closure in $\upsilon_{\mathcal{B}}(X)$ of an element in $\mathcal{B}$ is compact.
\begin{proof}
The proof of the first statement can be found in \cite{Gilman}. Furthermore, it is a well-known fact that the relatively pseudocompact subsets of $X$ are relatively compact in $\upsilon(X)$. The proof of the second statement is completely analogous.
\end{proof}
\end{prop}

\begin{defn}
A bornological universe $(X,\mathcal{B})$ will be called \emph{realcompact} if $\upsilon_{\mathcal{B}}(X)$ is equal to $X$.
\end{defn}

\begin{prop}
$(X,\mathcal{B})$ is realcompact iff $X$ is realcompact as a topological space and $\mathcal{B}$ contains only relatively compact subsets.
\begin{proof}
That this condition is necessary follows from the previous proposition. If this condition is satisfied, then each real valued, continuous map is automatically bounded on all elements of $\mathcal{B}$. This means that $\mathcal{C}_{\mathcal{B}}(X)$ is equal to $\mathcal{C}(X)$ and that $\upsilon_{\mathcal{B}}(X)$ is by definition equal to the Hewitt realcompactification of $X$. Since we assumed $X$ to be realcompact we obtain that $\upsilon_{\mathcal{B}}(X)$ equals $X$.
\end{proof}
\end{prop}

\begin{prop}
$(\upsilon_{\mathcal{B}}(X),\mathcal{C}(\upsilon_{\mathcal{B}}(X))$ is the completion of $(X,\mathcal{C}_{\mathcal{B}}(X))$.
\begin{proof}
Because $\upsilon_{\mathcal{B}}(X)$ is a realcompact Tychonoff space, we obtain that the uniform space $(\upsilon_{\mathcal{B}}(X),\mathcal{C}(\upsilon_{\mathcal{B}}(X))$ is complete (see \cite{Gilman}). To prove that it is the completion of $(X,\mathcal{C}_{\mathcal{B}}(X))$ we need to show that $(X,\mathcal{C}_{\mathcal{B}}(X))\hookrightarrow(\upsilon_{\mathcal{B}}(X),\mathcal{C}(\upsilon_{\mathcal{B}}(X))$ is a dense embedding. Because all $B\in\mathcal{B}$ are relatively compact in $\upsilon_{\mathcal{B}}(X)$ we have that the restriction of a map $f\in\mathcal{C}(\upsilon_{\mathcal{B}}(X))$ to $X$ is bounded and because the source $\mathcal{C}(\upsilon_{\mathcal{B}}(X))$ is initial this yields that $(X,\mathcal{C}_{\mathcal{B}}(X))\hookrightarrow(\upsilon_{\mathcal{B}}(X),\mathcal{C}(\upsilon_{\mathcal{B}}(X))$ is uniformly continuous. Moreover, because each $f\in\mathcal{C}_{\mathcal{B}}(X)$ can be extended to a real valued, continuous map on $\upsilon_{\mathcal{B}}(X)$ --- the extension being the restriction of $f^{\beta}$ to $\upsilon_{\mathcal{B}}(X)$ --- we obtain that this map is actually a dense embedding.
\end{proof}
\end{prop}

\begin{cor}
$(X,\mathcal{B})$ is realcompact iff $(X,\mathcal{C}_{\mathcal{B}}(X))$ is complete.
\end{cor}

By a \emph{character} on an algebra we mean a non-zero, real valued algebra morphism. One of the possible characterizations of realcompactness of a Tychonoff space $X$ is the following: for each character $\tau$ on $\mathcal{C}(X)$ there is an $x\in X$ such that $\tau(f)=f(x)$ for each $f\in \mathcal{C}(X)$. This means that for a realcompact topological space each character on the algebra of real valued, continuous maps is equal to a point-evaluation. In the setting for bornological universes we have the following result:

\begin{prop}
$(X,\mathcal{B})$ is realcompact iff each character on $\mathcal{C}_{\mathcal{B}}(X)$ is equal to a point-evaluation.
\begin{proof}
Suppose that this condition is satisfied and that $y$ is an element of $\upsilon_{\mathcal{B}}(X)$. This element $y$ defines a character on $\mathcal{C}_{\mathcal{B}}(X)$ by sending an $f$ to $f^{\beta}(y)$ and thus we can find an $x\in X$ such that $f(x)=f^{\beta}(y)$ for all $f\in\mathcal{C}_{\mathcal{B}}(X)$. Because the continuous maps from $\beta(X)$ into $[0,1]$ separate points we obtain that $x=y$.

If $X$ is equal to $\upsilon_{\mathcal{B}}(X)$, then all elements of $\mathcal{B}$ are relatively compact and therefore $\mathcal{C}_{\mathcal{B}}(X)$ is equal to $\mathcal{C}(X)$. If we combine this result with the fact that $\upsilon_{\mathcal{B}}(X)$, and therefore $X$, is realcompact as a topological space, then we get that each character on the algebra of bounded continuous maps on $(X,\mathcal{B})$ is equal to a point-evaluation.
\end{proof}
\end{prop}

For a Tychonoff space $X$ both compactness and realcompactness can be characterized using $z$-ultrafilters. In particular, $X$ is compact iff each $z$-ultrafilter is fixed and $X$ is realcompact iff each $z$-ultrafilter with the countable intersection property is fixed. We will see that for bornological universes $(X,\mathcal{B})$ with a normal underlying topological space $X$ the notion of realcompactness can be described with $z$-ultrafilters as well.

\begin{defn}
Let $(X,\mathcal{B})$ be a bornological universe. A decreasing sequence $(Z_n)_n$ of zero-sets is called \emph{unbounded} if for each $B\in\mathcal{B}$ there is an $n\in\mathbb{N}$ for which $\overline{B}\cap Z_n=\phi$.
\end{defn}

\begin{lem}\label{unboundedzultrafilter}
If $B$ is a subset of a normal space and $Z$ is a zero-set that intersects with all zero-sets that contain $B$, then $Z$ intersects with $\overline{B}$.
\begin{proof}
If this were not the case, then we could find a continuous map from $X$ into $[0,1]$ such that $f(B)\subseteq \{0\}$ and $f(Z)\subseteq \{1\}$. This would imply that there is a zero-set that contains $B$ and does not intersect with $Z$.
\end{proof}
\end{lem}

\begin{lem}
A subset $B$ of a topological space $X$ is relatively compact iff each open cover of $X$ contains a finite subcover of $B$.
\end{lem}

\begin{prop}
Let $(X,\mathcal{B})$ be a bornological universe with a normal underlying topological space $X$. $(X,\mathcal{B})$ is realcompact iff each $z$-ultrafilter that does not contain a decreasing, unbounded sequence of zero-sets, is fixed.
\begin{proof}
Let $(X,\mathcal{B})$ be realcompact and $F$ a $z$-ultrafilter on $X$ that does not contain a decreasing, unbounded sequence of zero-sets. We want to prove that $F$ has the countable intersection property. Take a decreasing sequence $(Z_n)_n$ of zero-sets in $F$. By definition we can find a $B\in\mathcal{B}$ such that for each $n\in\mathbb{N}$ we have  $\overline{B}\cap Z_n\not=\phi$. Because $\overline{B}$ is compact, the $z$-filter generated by the sets $\overline{B}\cap Z_n$ is fixed. Hence we obtain that $\bigcap_n Z_n$ is not empty. This means that $F$ has the countable intersection property and that, because $X$ is realcompact, $F$ is fixed.

To prove that this condition is sufficient we first take a $z$-ultrafilter $F$ on $X$ with the countable intersection property. Since each singleton in $X$ is bounded we obtain that $F$ does not contain any decreasing, unbounded sequences of zero-sets and thus, that $F$ is fixed. This yields that $X$ is realcompact. Now we take a $B\in\mathcal{B}$ and an open cover $\mathcal{G}$ of $X$. Let $\mathcal{Z}$ be the set of all zero-sets that contain a set $X\setminus G$ with $G\in\mathcal{G}$ and all zero-sets that contain $B$. If we assume that $\mathcal{G}$ does not contain a finite subcover of $B$, then $\mathcal{Z}$ has the finite intersection property and we can find a $z$-ultrafilter $F$ that contains $\mathcal{Z}$. From lemma \ref{unboundedzultrafilter} we obtain that each element of $F$ intersects with $\overline{B}$ and this of course implies that $F$ does not contain any decreasing, unbounded sequences of zero-sets. Because $F$ is not fixed we have to conclude that our original assumption was false and that $\mathcal{G}$ does contain a finite subcover of $B$.
\end{proof}
\end{prop}

When we look at the proof it is clear that if we do not assume the space $X$ to be normal, the condition stated in the previous proposition is still necessary.\\

The last results that we will encounter all concern Hausdorff uniformizable bornological universes $(X,\mathcal{B})$, i.e. $X$ is Tychonoff and $\mathcal{B}$ is saturated. We already saw that in a Hausdorff uniformizable bornological universe $(X,\mathcal{B})$ all relatively pseudocompact sets, and therefore all relatively compact sets, are contained in $\mathcal{B}$. We also know that in a realcompact bornological universe all bounded sets are relatively compact. This means that if a bornological universe $(X,\mathcal{B})$ is uniformizable and realcompact, the space $X$ is realcompact and $\mathcal{B}$ is the bornology of relatively compact sets (which is equal to the bornology of relatively pseudocompact sets for realcompact spaces). Since in this case all continuous, real valued maps on $X$ are bounded we automatically obtain that the converse implication is true as well.

\begin{prop}
$(X,\mathcal{B})$ is Hausdorff uniformizable and realcompact iff it is isomorphic to a closed subspace of a product of real lines.
\begin{proof}
Each realcompact topological space $X$ is isomorphic to a closed subset $C$ of a product $\mathbb{R}^{\alpha}$ of real lines. In such a space a subset $B$ is bounded iff each of its projections is bounded or, equivalently, it is contained in a compact subset of $\mathbb{R}^{\alpha}$. Since $C$ is closed this is equivalent to saying that the closure of $B$ in $C$ is compact. This means that this condition is necessary.

We know that a topological space that is isomorphic to a closed subset of a product of real lines is realcompact and that the closure of a bounded set in such a space is compact, so this condition is also sufficient.
\end{proof}
\end{prop}

\begin{defn}\label{realcompactification}
For a Hausdorff uniformizable bornological universe $(X,\mathcal{B})$ we define $\upsilon(X,\mathcal{B})$ as the object $(\upsilon_{\mathcal{B}},\upsilon(\mathcal{B}))$, where $\upsilon(\mathcal{B})$ denotes the set of all relatively compact subsets of $\upsilon_{\mathcal{B}}(X)$. This is again a Hausdorff uniformizable bornological universe and by definition it is also realcompact.
\end{defn}

\begin{prop}
Let $f:(X,\mathcal{B})\to (Y,\mathcal{B}')$ be a morphism. If we define $\upsilon(f)$ as the restriction of $\beta(f)$ to $\upsilon_{\mathcal{B}}(X)$, then $\upsilon(f):\upsilon(X,\mathcal{B})\to \upsilon(Y,\mathcal{B}')$ is a morphism.
\begin{proof}
What we need to prove is that $\beta(f)$ is a morphism from $\upsilon_{\mathcal{B}}(X)$ to $\upsilon_{\mathcal{B}'}(Y)$. Let $x$ be an element of $\upsilon_{\mathcal{B}}(X)$. For each $g\in\mathcal{C}_{\mathcal{B}'}(Y)$ we have that $g\circ f$ is an element of $\mathcal{C}_{\mathcal{B}}(X)$ and thus that $(g\circ f)^{\beta}(x)$ is real. Because the map $(g\circ f)^{\beta}$ is equal to $g^{\beta}\circ \beta(f)$ we obtain that $\beta(f)(x)$ is in $\upsilon_{\mathcal{B}'}(Y)$. That $\upsilon(f)$ is bounded is trivially true, since relatively compact sets are mapped to relatively compact sets. 
\end{proof}
\end{prop}

\begin{cor}
$\upsilon(X,\mathcal{B})$ is the realcompact reflection of $(X,\mathcal{B})$ in the category of Hausdorff uniformizable bornological universes and bounded, continuous maps.
\begin{proof}
We know that $(X,\mathcal{B})\hookrightarrow \upsilon(X,\mathcal{B})$ is a morphism. Furthermore, if $(Y,\mathcal{B}')$ is realcompact, then $\upsilon(Y,\mathcal{B})$ is equal to $(Y,\mathcal{B})$ and thus we obtain for each map $f:(X,\mathcal{B})\to (Y,\mathcal{B}')$ that is bounded and continuous that $\upsilon(f)$ is a morphism from $\upsilon(X,\mathcal{B})$ to $(Y,\mathcal{B})$.
\end{proof}
\end{cor}

The following proposition extends the well-known theorem that states that the Hewitt realcompactification of a dense subspace $Y$ of $X$ equals the Hewitt realcompactification of $X$ iff $Y$ is $C$-embedded in $X$.

\begin{prop}
Let $Y$ be a dense subspace of $X$ and let $\mathcal{B}'$ be the subspace bornology on $Y$ derived from $\mathcal{B}$. The following statements are equivalent:
\begin{enumerate}
\item if $(Z,\mathcal{B}'')$ is realcompact, then each map $f:(Y,\mathcal{B}')\to (Z,\mathcal{B}'')$ that is continuous and bounded has an extension to $(X,\mathcal{B})$,
\item each realvalued map $f$ on $(Y,\mathcal{B}')$ that is continuous and bounded has an extension to $(X,\mathcal{B})$,
\item $\upsilon(Y,\mathcal{B}')=\upsilon(X,\mathcal{B})$.
\end{enumerate}
\begin{proof}
The first statement implies the second because $\mathbb{R}$ endowed with the bornology of relatively compact subsets is a realcompact bornological universe.

From the second statement we obtain that each realvalued, continuous map $f$ that is bounded on $Y$ has an extension to $X$ . This implies that $\beta(Y)=\beta(X)$. By definition we have that $\upsilon_{\mathcal{B}'}(Y)\subseteq \upsilon_{\mathcal{B}}(X)$. Now take an $x\in\upsilon_{\mathcal{B}}(X)$, an $f'\in\mathcal{C}_{\mathcal{B}'}(Y)$ and let $f$ be the extension of $f'$ to $X$. Since they coincide on a dense subset we have that $(f')^{\beta}$ equals $f^{\beta}$. Hence we get that $(f')^{\beta}(x)$ is equal to $f^{\beta}(x)$ and therefore an element of the reals. Because this is true for all maps in $\mathcal{C}_{\mathcal{B}'}(Y)$ we obtain that $\upsilon_{\mathcal{B}}(X)$ is a subset of $\upsilon_{\mathcal{B}'}(Y)$.  From the foregoing we obtain that $\upsilon_{\mathcal{B}}(X)$ is equal to $\upsilon_{\mathcal{B}'}(Y)$ and thus that $\upsilon(X,\mathcal{B})$ equals $\upsilon(Y,\mathcal{B}')$.

Suppose the third statement is true and take a map $f:(Y,\mathcal{B}')\to (Z,\mathcal{B}'')$, where $f$ is bounded and continuous and $(Z,\mathcal{B}'')$ is realcompact. We know that the map $\upsilon(f):\upsilon(Y,\mathcal{B})\to \upsilon(Z,\mathcal{B}'')$ is bounded and continuous and that $\upsilon(Z,\mathcal{B}'')$ is equal to $(Z,\mathcal{B}'')$. This implies that the restriction of $\upsilon(f)$ to $X$ is a map to $(Z,\mathcal{B}'')$ that extends $f$.
\end{proof}
\end{prop}


\begin{thebibliography}{99}
\bibitem{Beer1} G. Beer. Embeddings of bornological universes. \textit{Set-valued Analysis 16}. 477-488.  Springer, Dordrecht, (2008).
\bibitem{Beer5} G. Beer, S. Naimpally, J. Rodr\'iguez-L\'opez. $\mathcal{S}$-topologies and bounded convergences. \textit{Journal of Mathematical Analysis and Applications 339}. 542-552. Elsevier, San Diego, (2008).
\bibitem{Beer6} G. Beer, S. Levi. Gap, excess and bornological convergence. \textit{Set-Valued Analysis 16}. 489-506. Springer, Dordrecht, (2008).
\bibitem{Beer7} G. Beer, S. Levi. Pseudometrizable bornological convergence is Attouch-Wets convergence. \textit{Journal of Convex Analysis 15}. 439-453. Heldermann Verlag, Lemgo, (2008).
\bibitem{Beer2} G. Beer, S. Levi. Total boundedness and bornologies. \textit{Topology and its applications 156}. 1271-1288. Elsevier, Amsterdam, (2009).
\bibitem{Beer3} G. Beer, M. Segura. Well-posedness, bornologies, and the structure of metric spaces. \textit{Applied general topology 10}. 131-157. Universidad PolitŽcnica de Valencia, Valencia, (2009).
\bibitem{Beer4} G. Beer, S. Levi. Strong uniform continuity. \textit{Journal of Mathematical Analysis and Applications 350}. 568-589. Elsevier, San Diego, (2009).
\bibitem{Bourbaki} N. Bourbaki. \textit{Topologie G\'en\'erale}. Hermann, Paris, (1965).
\bibitem{Gilman} L. Gilman, M. Jerison. \textit{Rings of Continuous Functions}. D. Van Nostrand Company, New York, (1960).
\bibitem{Hejcman}  J. Hejcman. Boundedness in uniform spaces and topological groups. \textit{Czechoslovak Mathematical Journal 84}. 544-563. Springer, Dordrecht, (1959).
\bibitem{Hu1} S. Hu. Boundedness in a topological space. \textit{Journal de Math\'ematiques Pures et Appliqu\'ees 28}. 287-320. Elsevier, Amsterdam, (1949).
\bibitem{Hu2} S. Hu. \textit{Introduction to General Topology}. Holden-Day, San Francisco, (1966).
\bibitem{Schmets} J. Schmets. Espaces de fonctions continues. \textit{Lecture Notes in Mathematics 519}. Springer-Verlag, Berlin, (1976).
\bibitem{Lechicki} A. Lechicki, S. Levi, A. Spakowski. Bornological convergences. \textit{Journal of Mathematical Analysis and Applications 297}. 751-770. Elsevier, San Diego, (2004).
\bibitem{Vroegrijk1} T. Vroegrijk. Pointwise bornological spaces. \textit{Topology and its Applications 156}. 2019-2027. Elsevier, Amsterdam, (2009).
\bibitem{Vroegrijk2} T. Vroegrijk. Pointwise bornological vector spaces. \textit{Topology and its applications}. to appear.
\end{thebibliography}
\end{document}